\documentclass[12pt,paper]{amsart}
\usepackage{amsmath}
\usepackage{amsxtra}
\usepackage{amscd}
\usepackage{amsthm}
\usepackage{amsfonts}
\usepackage{amssymb}
\usepackage{eucal}
\usepackage{mathrsfs}
\usepackage{hyperref}

\def\P{{\mathbb P}}

\def\Z{{\mathbb Z}}

\newtheorem{theorem}{Theorem}[section]
\newtheorem{definition}{Definition}[section]
\newtheorem{lemma}[theorem]{Lemma}

\newtheorem{proposition}[theorem]{Proposition}
\newtheorem{corollary}[theorem]{Corollary}

\title{Some effects of Veronese map on syzygies of projective varieties}
\author[E. Park]{Euisung Park}
\address {Euisung Park : School of Mathematics, Korea Institute for Advanced Study,
207-43 Cheongryangri-dong, Dongdaemun-gu, Seoul 130-722, Republic
of Korea,} \email{puserdos@kias.re.kr}
\thanks{Mathematics Subject Classification (2000): 13D02, 14N15}

\begin{document}

\thispagestyle{empty} \maketitle

\begin{abstract}
Let $X \subset \P^r$ be a nondegenerate projective variety and let
$\nu_{\ell} : \P^r \rightarrow \P^N$ be the $\ell$-th Veronese
embedding. In this paper we study the higher normality, defining
equations and syzygies among them for the projective embedding
$\nu_{\ell} (X) \subset \P^N$. We obtain that for a very ample
line bundle $L \in \mbox{Pic}X$ such that $X \subset \P H^0 (X,L)$
is $m$-regular in the sense of Castelnuovo-Mumford, $(X,L^{\ell})$
satisfies property $N_{\ell}$ for all $\ell \geq m$ (Theorem
\ref{thm:main1}). This is a generalization of M. Green's work that
$(\P^r,\mathcal{O}_{\P^r} (\ell))$ satisfies property $N_{\ell}$.
Also our result refines works of Ein-Lazarsfeld in \cite{EL},
Gallego-Purnaprajna in \cite{GP} and E. Rubei in \cite{Rubei1},
\cite{Rubei} and \cite{Ru1}.
\end{abstract}

\tableofcontents \setcounter{page}{1}

\section{Introduction}
Let $X \subset \P^r$ be a nondegenerate projective variety where
the embedding is defined by an $(r+1)$-dimensional subspace $V
\subset H^0 (X,\mathcal{O}_X (1))$ and hence $\P^r = \P (V)$. Let
$V_{\ell}$ be the image of the natural homomorphism
\begin{equation*}
S^{\ell} V \cong H^0 (\P^r,\mathcal{O}_{\P^r} (\ell)) \rightarrow
H^0 (X,\mathcal{O}_X (\ell)).
\end{equation*}
In this paper we are concerned with the problem to see which
information on
\begin{equation*}
X \subset \P (V_{\ell})~\mbox{and}~X \subset \P H^0
(X,\mathcal{O}_X (\ell))
\end{equation*}
can be deduced from $X \subset \P (V)$. The first result to this
problem is due to Mumford.

\begin{theorem}[D. Mumford, \cite{Mumford}]\label{thm:Mumford}
Let $d$ be the degree of a nondegenerate projective variety $X
\subset \P^r$. If $\ell \geq d$, then $X \subset \P (V_{\ell})$
and $X \subset \P H^0 (X,\mathcal{O}_X (\ell))$ are represented as
an intersection of quadratic forms.
\end{theorem}

\noindent Later the result about the embedding defined by $H^0
(X,\mathcal{O}_X (\ell))$ is generalized to more general
statements about higher syzygies. For precise statements, we
recall the definition of property $N_p$.

\begin{definition}
Let $X$ be a projective variety and let $L \in \mbox{Pic}X$ a very
ample line bundle defining an embedding $X \hookrightarrow \P H^0
(X,L)$. Denote by $S=Sym^\bullet H^0 (X,L)$ the homogeneous
coordinate ring of $\P H^0 (X,L)$, and consider the graded
$S$-module $R(L)=\oplus_{n \in \Z} H^0 (X,L^n)$. Let
\begin{eqnarray*}
\cdots \rightarrow \oplus_j S^{\beta_{i,j}}(-i-j) \rightarrow
\cdots \rightarrow \oplus_j S^{\beta_{1,j}}(-1-j) \rightarrow
\oplus_j S^{\beta_{0,j}}(-j) \rightarrow R(L) \rightarrow 0
\end{eqnarray*}
be a minimal graded free resolution of $R(L)$. The line bundle $L$
is said to satisfy \textit{property} $N_p$ if $\beta_{i,j}=0$ for
$0 \leq i \leq p$ and $j \geq 2$.
\end{definition}

\noindent Therefore property $N_0$ holds if and only if $X
\hookrightarrow \P H^0 (X,L)$ is a projectively normal embedding,
property $N_1$ holds if and only if property $N_0$ is satisfied
and the homogeneous ideal is generated by quadrics, and property
$N_p$ holds for $p \geq 2$ if and only if it has property $N_1$
and the $k^{th}$ syzygies among the quadrics are generated by
linear syzygies for all $1 \leq k \leq p-1$. For
$(X,L)=(\P^n,\mathcal{O}_{\P^n} (1))$, M. Green obtain

\begin{theorem}[M. Green, \cite{Green2}]\label{thm:GreenVeronese}
$(\P^n,\mathcal{O}_{\P^n} (\ell))$ satisfies property $N_{\ell}$.
\end{theorem}

\noindent Also L. Ein and R. Lazarsfeld generalize Theorem
\ref{thm:Mumford} as follows:

\begin{theorem}[L. Ein and R. Lazarsfeld, \cite{EL}]\label{thm:EL}
Let $X$ be a smooth complex projective variety and let $L \in
\mbox{Pic}X$ be a very ample line bundle of degree $d$. Then
$(X,L^{\ell})$ satisfies property $N_{\ell+1-d}$ for all $\ell
\geq d-1$.
\end{theorem}

\noindent Later F. J. Gallego and B. P. Purnaprajna obtain the
following result:

\begin{theorem}[F. J. Gallego and B. P. Purnaprajna, \cite{GP}]\label{thm:GP}
Let $X$ be a projective variety and let $L \in \mbox{Pic}X$ be a
base point free and ample line bundle. Assume that $\mathcal{O}_X$
is $m$-regular with respect to $L$, i.e., $H^i (X,L^{m-i})=0$ for
all $i \geq 1$. Then for $\ell \geq \mbox{max} \{m-1+p,m+1,p+1\}$,
$(X,L^{\ell})$ satisfies property $N_p$. In particular, if $m \geq
2$ and $\ell \geq 2m-1$, then $(X,L^{\ell})$ satisfies property
$N_{\ell-m+1}$.
\end{theorem}

\noindent {\bf Remark 1.1.} Theorem \ref{thm:GP} extends Theorem
\ref{thm:EL}. Indeed let $L$ be a very ample line bundle of degree
$d$ on a smooth complex projective variety $X$ of dimension $n$.
Then the line bundle $D=L^{d-n-2} \otimes K_X ^{-1}$ is base point
free. One can find the details at the proof of Proposition 3.3 in
\cite{EL}. Thus $L^{\ell} = K_X \otimes L^{\ell+n+2-d} \otimes D$
and hence $\mathcal{O}_X$ is $(d-1)$-regular with respect to $L$
by Kodaira vanishing theorem. If $d \geq 2$, then Theorem
\ref{thm:GP} implies that\\

$(X,L^{\ell})$ satisfies $\begin{cases} \mbox{property
$N_{\ell-1}$} &
\mbox{for $d=2$ and $\ell \geq 2$, and} \\
\mbox{property $N_{\ell+2-d}$} & \mbox{for $d \geq 3$ and $\ell
\geq d-1$.} \end{cases}$\\

\noindent If $d=1$, then $(X,L) = (\P^n,\mathcal{O}_{\P^n} (1))$.
In this case, see Theorem \ref{thm:GreenVeronese}. \qed \\

Roughly speaking, these works guarantee that for a very ample line
bundle $L$ on a projective variety $X$, $(X,L^{\ell})$ satisfies
property $N_{\ell -c}$ for all sufficiently large $\ell$ where $c
\geq 0$ is an invariant of $(X,L)$. And our first main theorem
shows that $(X,L^{\ell})$ satisfies property $N_{\ell}$ for all
sufficiently large $\ell$.

\begin{theorem}\label{thm:main1}
Let $X$ be a projective variety and let $L \in \mbox{Pic}X$ be a
very ample line bundle such that $X \subset \P H^0 (X,L)=\P$ is
$m$-regular in the sense of Castelnuovo-Mumford, that is, $H^i
(\P,\mathcal{I}_{X/\P} (m-i))=0$ for all $i \geq 1$.
Then \\
$(1)$ $L^{\ell}$ is normally generated for $\ell \geq
\frac{m-1}{2}$.\\
$(2)$ $L^{m-1}$ satisfies property $N_{m-2}$.\\
$(3)$ $L^{\ell}$ satisfies property $N_{\ell}$ for all $\ell \geq
m$.
\end{theorem}

This generalizes Green's theorem for $(\P^n,\mathcal{O}_{\P^n}
(1))$ to arbitrary $(X,L)$. Also this refines Theorem \ref{thm:EL}
and \ref{thm:GP}. In $\S 4$ we apply our theorems to various
projective varieties. We refine or reprove some known facts about
property $N_p$ of K3 surfaces, Enriques surfaces and abelian
varieties in \cite{GP}, \cite{Rubei1} and \cite{Ru1}.\\

Theorem \ref{thm:main1} is a consequence of a more general result
on some effects of the Veronese map on higher normality, defining
equations and syzygies among them. Let $X \subset \P (V)$ be a
nondegenerate projective variety defined by a subspace $V \subset
H^0 (X,\mathcal{O}_X (1))$ and let $V_{\ell}$ be the image of the
homomorphism
\begin{equation*}
S^{\ell} V \cong H^0 (\P^r,\mathcal{O}_{\P^r} (\ell)) \rightarrow
H^0 (X,\mathcal{O}_X (\ell)).
\end{equation*}
From the surjective map $H^0 (\P^r,\mathcal{O}_{\P^r} (\ell))
\rightarrow V_{\ell}$, we can consider $\P (V_{\ell})$ as a linear
subspace of $\P (S^{\ell}V)$. Thus we get the following
commutative diagram:
\begin{equation*}
\begin{CD}
\quad X              & \quad  \subset \quad & \P (V) \\
\nu_{\ell} \downarrow     &                      & \downarrow \nu_{\ell} \\
\P (V_{\ell})  & \quad  \subset \quad & \P (S^{\ell} V)
\end{CD}
\end{equation*}
We call $X \subset \P (V_{\ell})$ the $\ell$-th Veronese
embedding. Recall that as defined in \cite{EGHP}, a projective
variety embedded in a projective space satisfies property
$N_{2,1}$ if the homogeneous ideal is generated by quadrics and it
satisfies property $N_{2,p}$ for $p \geq 2$ if property $N_{2,1}$
holds and the $k$-th syzygies among the quadrics are generated by
linear syzygies for all $1 \leq k \leq p-1$. See $\S 2.3$. We
prove the following

\begin{theorem}\label{thm:veronese}
Let $X \subset \P (V)$ be a nondegenerate projective variety
defined by a subspace $V \subset H^0 (X,\mathcal{O}_X (1))$. Let
$s,t$ and $m$ be integers such that
\begin{enumerate}
\item[$(\alpha)$] $k$-normality holds for all $\ell \geq s$,
\item[$(\beta)$] the homogeneous ideal is generated by forms of
degree $\leq t$, and \item[$(\gamma)$] $X \subset \P (V)$ is
$m$-regular, i.e., $H^1 (\P (V),\mathcal{I}_X (m-i))=0$ for all $i
\geq 1$.
\end{enumerate}
Then for the $\ell$-th Veronese embedding, \\
$(1)$ $\nu_{\ell} (X)$ satisfies $k$-normality for all $k \geq \frac{s}{\ell}$.\\
$(2)$ The homogeneous ideal of $\nu_{\ell} (X)$ is generated by
forms of degree $\leq \mbox{max} \{2,\frac{t}{\ell}\}$.

In particular, $\nu_{\ell} (X)$ is ideal-theoretically cut out
by quadrics if $\ell \geq \frac{t}{2}$.\\
$(3)$ $\nu_{\ell} (X)$ satisfies $\begin{cases} \mbox{property
$N_{2,2\ell-m}$} &  \mbox{for $\frac{m+1}{2} \leq \ell < m$, and}\\
\mbox{property $N_{2,\ell}$} &  \mbox{for $\ell \geq m$.}
\end{cases}$
\end{theorem}

For the proof we first investigate syzygies among the defining
equations of the ``degenerate" embedding $X \subset \P
(S^{\ell}V)$. Let $\mathcal{I}_{X/\P (V)}$, $\mathcal{I}_{\P (V) /
\P (S^{\ell}V)}$, $\mathcal{I}_{X/ \P (S^{\ell}V)}$ and
$\mathcal{I}_{X/ \P (V_\ell)}$ be sheaves of ideals of $X \subset
\P (V)$, $\nu_{\ell} (\P (V)) \subset \P (S^{\ell}V)$,
$\nu_{\ell}(X) \subset \P (S^{\ell}V)$ and $\nu_{\ell}(X) \subset
\P (V_{\ell})$, respectively. Then we have the short exact
sequence
\begin{equation*}
0 \rightarrow  \mathcal{I}_{\P (V) / \P (S^{\ell}V)} \rightarrow
\mathcal{I}_{X/ \P (S^{\ell}V)} \rightarrow  \mathcal{I}_{X/\P
(V)} \rightarrow 0.
\end{equation*}
For details, see the proof of Theorem \ref{thm:veronese}. For
$\mathcal{I}_{\P (V) / \P (S^{\ell}V)}$ we know that $\nu_{\ell}
(\P (V)) \subset \P (S^{\ell}V)$ satisfies property $N_{\ell}$.
Also our assumptions $(\alpha)\sim(\gamma)$ give some information
about $\mathcal{I}_{X/\P (V)}$. Thus we obtain some results about
syzygy modules of the degenerate embedding $X \subset \P (S^{\ell}
V)$. As discussed in $\S 2.3$, syzygies among defining equations
of $X \subset \P (S^{\ell}V)$ are closely related to those of $X
\subset \P (V_{\ell})$. These observations enables us to prove
Theorem \ref{thm:veronese}.\\

The proof of Theorem \ref{thm:veronese}.$(3)$ combines the theory
of Castelnuovo-Mumford regularity with Green's theorem. Indeed we
obtain the following more general statement:

\begin{theorem}\label{thm:imply}
Let $L$ be a very ample line bundle on a projective variety $X$
such that $X \subset \P H^0 (X,L)=\P^r$ is $m$-regular. For $\ell
\geq m$, suppose that $(\P^r,\mathcal{O}_{\P^r}(\ell))$ satisfies
property $N_p$. Then $(X,L^{\ell})$ satisfies property $N_p$.
\end{theorem}

This result guarantees that if Green's theorem is refined, then
Theorem \ref{thm:main1} is automatically refined. See the
following two examples.\\

\noindent {\bf Example 1.1.} In \cite{Ru2}, E. Rubei shows that
$(\P^r,\mathcal{O}_{\P^r} (3))$ satisfies property $N_4$.
Therefore if $X \subset \P H^0 (X,L)$ is $3$-regular, then
$(X,L^3)$ satisfies property $N_4$. See $\S 4.3$ for applications
to $3$-regular varieties.   \qed \\

\noindent {\bf Example 1.2.} Let $X \subset \P^r$ be a quadratic
hypersurface or a rational normal scroll. Then $(X,\mathcal{O}_X
(2))$ satisfies property $N_5$. Indeed T. Jozefiak, P. Pragacz and
J. Weyman prove that $(\P^n,\mathcal{O}_{\P^n} (2))$, $n \geq 3$,
satisfies property $N_p$ if and only if $p \leq 5$\cite{JPW}.
Since $X \subset \P^r$ is $2$-regular, the assertion comes from
Theorem \ref{thm:imply}. \qed \\

\noindent {\bf Example 1.3.} Let $C \subset \P^2$ be a plane cubic
curve which is not necessarily smooth or irreducible. When $C$ is
smooth, $(C,\mathcal{O}_C (\ell))$ satisfies property
$N_{3\ell-3}$ by Green's ``$2g+1+p$" theorem in \cite{Green}. By
Theorem \ref{thm:imply}, we can extend this fact for arbitrary
cubic plane curves. Indeed $(\P^2,\mathcal{O}_{\P^2} (\ell))$,
$\ell \geq 3$, satisfies property $N_{3\ell-3}$. For details, we
refer the reader to see \cite{OP}. Therefore $(C,\mathcal{O}_C
(\ell))$ satisfies property $N_{3\ell-3}$, $\ell \geq 3$,
satisfies property $N_{3\ell-3}$.  \qed \\

Recently G. Ottaviani and R. Paoletti investigate syzygies of
Veronese embedding of projective spaces. By Ottaviani-Paoletti's
theorem in \cite{OP}, if $r \geq 2$ and $\ell \geq 3$ then
$(\P^r,\mathcal{O}_{\P^r} (\ell))$ fails to satisfy property
$N_{3\ell-2}$.
Also they formulate the following\\

{\bf Conjecture A.} Assume that $r \geq 2$ and $\ell \geq 3$. Then
$(\P^r,\mathcal{O}_{\P^r} (\ell))$ satisfies property $N_p$
if and only if $p \leq 3\ell-3$.\\

\noindent Ottaviani-Paoletti's conjecture can be generalized as follows:\\

{\bf Conjecture B.} Let $L$ be a very ample line bundle on a
projective variety $X$ such that $X \subset \P H^0 (X,L)$ is
$m$-regular. Then for $\ell \geq m$, $(X,L^{\ell})$
satisfies property $N_{3\ell -3}$.\\

\noindent Clearly Theorem \ref{thm:imply} guarantees that if
``Conjecture A" is true, then ``Conjecture B" is also true. This
conjecture says that for all sufficiently large $\ell$, syzygies
of $(X,L^{\ell})$ is closed related to those of
$(\P^n,\mathcal{O}_{\P^n}(\ell))$. Also the following example
shows that the assumption ``$\ell \geq m$" is necessary.  \\

\noindent {\bf Example 1.4.} Let $X$ be an elliptic ruled surface
with the numerical invariant $e \geq 0$. Then $aC_0+bf$ is very
ample if and only if $a \geq 1$ and $b-ae \geq 3$. Also $aC_0+bf$
satisfies property $N_p$ if and only if $a \geq 1$ and $b-ae \geq
3+p$ by Theorem 1.4 in \cite{ES}. Therefore Conjecture B is true
for every very ample line bundle on $X$. Now let $L=aC_0 +bf$ be
such that $a\geq 1$ and $b-ae=3$. Then $L$ is very ample and $X
\subset \P H^0 (X,L)$ is $3$-regular. Note that $(X,L^2)$
satisfies property $N_p$ if and only if $p \leq 3$ while
$(\P^n,\mathcal{O}_{\P^n}(2))$, $n\geq 3$, satisfies property
$N_p$ if and only if $p \leq 5$. \qed \\

{\bf Organization of the paper.} In $\S 2$, we review some basic
facts to study syzygies. In $\S 3$ we deal with higher normality
and defining equations of Veronese map of arbitrary projective
varieties. Also we prove our main results on syzygies of powers of
a very ample line bundle. Finally we present some
applications of our results in $\S 4$.\\

{\bf Notations.} Throughout this paper all varieties are defined
over an algebraically closed field $K$ with $\mbox{char}K =0$. For
a finite dimensional $K$-vector space $V$, $\P(V)$ denotes the
projective space of one-dimensional quotients of $V$.\\

{\bf Acknowledgement.} I thank Elena Rubei for useful comments
about $\S 4.5$.\\

\section{Preliminaries}
\subsection{Surjection of multiplication maps.} To prove our main
result we need to show that a particular multiplication map of
global sections on coherent sheaves is surjective. Here we present
an elementary but useful lemma for this aim.

\begin{lemma}\label{lem:mult}
Let $X$ be a projective variety and let $\mathcal{E}$ be a locally
free sheaf on $X$. For a short sequence $0 \rightarrow \mathcal{F}
\rightarrow \mathcal{G} \rightarrow \mathcal{H} \rightarrow 0$ of
coherent sheaves on $X$, assume that

$(\alpha)$ $\mathcal{F}$ and $\mathcal{H}$ are globally generated
(and hence $\mathcal{G}$ is globally generated),

$(\beta)$ the following two sequences are exact, and
$$\begin{cases}
0 \rightarrow H^0 (X,\mathcal{F}) \rightarrow H^0 (X,\mathcal{G})
\rightarrow H^0 (X,\mathcal{H}) \rightarrow 0  \\
0 \rightarrow H^0 (X,\mathcal{F} \otimes \mathcal{E}) \rightarrow
H^0 (X,\mathcal{G}\otimes \mathcal{E}) \rightarrow H^0
(X,\mathcal{H}\otimes \mathcal{E}) \rightarrow 0
\end{cases}$$

$(\gamma)$ the following two multiplication maps are surjective.
$$\begin{cases}
H^0 (X,\mathcal{F}) \otimes H^0 (X,\mathcal{E}) \rightarrow H^0
(X,\mathcal{F} \otimes \mathcal{E})  \\
H^0 (X,\mathcal{H}) \otimes H^0 (X,\mathcal{E})  \rightarrow H^0
(X,\mathcal{H} \otimes \mathcal{E})
\end{cases}$$
Then the multiplication map $H^0 (X,\mathcal{G}) \otimes H^0
(X,\mathcal{E}) \rightarrow H^0 (X,\mathcal{G} \otimes
\mathcal{E})$ is surjective.
\end{lemma}

\begin{proof}
By our assumption we have the commutative diagram
\begin{equation*}
\begin{CD}
0 \rightarrow & H^0 (X,\mathcal{F})\otimes H^0 (X,\mathcal{E}) & \rightarrow & H^0 (X,\mathcal{G}) \otimes H^0 (X,\mathcal{E}) & \rightarrow & H^0 (X,\mathcal{H})\otimes H^0 (X,\mathcal{E}) & \rightarrow 0 \\
              & \downarrow             &             &  \downarrow            &             &             \downarrow &         \\
0 \rightarrow & H^0 (X,\mathcal{F}\otimes\mathcal{E})  &
\rightarrow & H^0 (X,\mathcal{G}\otimes\mathcal{E}) & \rightarrow
& H^0 (X,\mathcal{H}\otimes\mathcal{E}) & \rightarrow 0.
\end{CD}
\end{equation*}
where the two rows are exact. The assertion comes from snake
lemma.
\end{proof}

\subsection{Cohomological criterion of property $N_{2,p}$} Let $X \subset \P^r$ be a
nondegenerate projective variety and let $I_X$ be the homogeneous
ideal of $X$ in the homogeneous coordinate ring $S$ of $\P^r$. For
$p \geq 1$, $X$ is said to satisfies property $N_{2,p}$ if $I_X$
admits the minimal free resolution of the form
\begin{equation*}
\cdots \rightarrow S(-p-1)^{\beta_{p-1,2}}  \rightarrow \cdots
\rightarrow  S(-3)^{\beta_{1,2}} \rightarrow S(-2)^{\beta_{0,2}}
\rightarrow I_X \rightarrow 0.
\end{equation*}
Therefore property $N_{2,1}$ holds if $I_X$ is generated by
quadrics. For $p \geq 2$, it is said to satisfy property $N_{2,p}$
if property $N_{2,1}$ holds and the $k$-th syzygies among the
quadrics are generated by linear syzygies for all $1 \leq k \leq
p-1$. Equivalently, property $N_{2,p}$ holds if $\mbox{Tor}^S _i
(I_X,K)$ is a vector space concentrated in degrees $\leq i+2$ for
all $i \leq p-1$. It is well known that $\mbox{Tor}^S _i (I_X,K)$
can be read off as follows:

\begin{theorem}\label{thm:exactsequence}
Put $V = H^0 (\P^r,\mathcal{O}_{\P^r} (1))$ and
$\mathcal{M}=\Omega _{\P^r} (1)$. Then for the sheaf of ideals
$\mathcal{I}_X$ of $X$, there is an exact sequence
\begin{equation*}
\begin{CD}
\bigwedge^{i+1}V \otimes H^0 ( \P^r, \mathcal{I}_X (j-1))
\stackrel{\alpha_{i,j}}{\rightarrow} H^0 (\P^r,\bigwedge^i
\mathcal{M}
\otimes \mathcal{I}_X (j) )  \rightarrow \\
H^1 (\P,\bigwedge^{i+1} \mathcal{M} \otimes \mathcal{I}_X (j-1) )
\rightarrow  \bigwedge^{i+1} V \otimes H^1 ( \P, \mathcal{I}_X
(j-1))
\end{CD}
\end{equation*}
with $\mbox{Coker}(\alpha_{i,j})=Tor^{i+j} _i (I_X,K)$. Therefore
we have the exact sequence
\begin{equation*}
0 \rightarrow Tor^{i+j} _i (I_X,K) \rightarrow  H^1 (\P^r
,\bigwedge^{i+1} \mathcal{M} \otimes \mathcal{I}_X (j-1) )
\rightarrow \bigwedge^{i+1} V \otimes H^1 ( \P^r,\mathcal{I}_X
(j-1)).
\end{equation*}
\end{theorem}

\begin{proof}
See Theorem (1.b.4) in \cite{Green} or Theorem 4.5 in
\cite{Eisenbud}.
\end{proof}

\begin{corollary}\label{cor:criterionN2p}
$(1)$ If $H^1 (\P^r ,\bigwedge^i \mathcal{M} \otimes \mathcal{I}_X
(j) )=0$ for $1 \leq i \leq p$ and $j \geq 2$, then $X \subset
\P^r$
satisfies property $N_{2,p}$.\\
$(2)$ If $X \subset \P^r$ is projectively normal, then it
satisfies property $N_{2,p}$ if and only if $H^1 (\P^r
,\bigwedge^i \mathcal{M} \otimes \mathcal{I}_X (j) )=0$ for $1
\leq i \leq p$ and $j \geq 2$.
\end{corollary}

\begin{proof}
By definition property $N_{2,p}$ holds if and only if
\begin{equation*}
Tor^{i+j} _i (I_X,K)=0\quad \mbox{for $0 \leq i \leq p-1$ and $j
\geq 3$.}
\end{equation*}
Thus the assertion comes immediately from the exact sequence in
Theorem \ref{thm:exactsequence}.
\end{proof}

\subsection{Syzygies of degenerate varieties}
Let $X \subset \P^{r+e}(e \geq 1)$ be a degenerate projective
variety, i.e., there exists a linear subspace $\Lambda \cong \P^r
\subset \P^{r+e}$ such that $X \subset \Lambda$. To simplify
notations, put $\mathcal{M}_n = \Omega_{\P^n} (1)$. It is easily
checked that
\begin{equation*}
\mathcal{M}_{r+e}|_{\Lambda} \cong \mathcal{M}_r \oplus
\mathcal{O}_{\P^r} ^e.
\end{equation*}
Let $\mathcal{I}_{X/\P^{r+e}}$, $\mathcal{I}_{\Lambda/\P^{r+e}}$
and $\mathcal{I}_{X/\P^r}$ be sheaves of ideals of $X \subset
\P^{r+e}$, $\Lambda \subset \P^{r+e}$ and $X \subset \P^r$,
respectively.

\begin{lemma}\label{lem:degenrate}
Under the situation just stated, suppose that for some $k \geq 0$
and $j \geq 1$,
\begin{equation*}
H^1 (\P^{r+e},\bigwedge^k \mathcal{M}_{r+e} \otimes
\mathcal{I}_{X/\P^{r+e}} (j)) =0
\end{equation*}
Then $H^1 (\P^r,\bigwedge^k \mathcal{M}_r \otimes
\mathcal{I}_{X/\P^r} (j)) =0$.
\end{lemma}

\begin{proof}
We may assume that $e=1$ and hence
$\mathcal{I}_{\Lambda/\P^{r+1}}=\mathcal{O}_{\P^{r+1}}(-1)$. By
snake lemma we have the following commutative diagram:
\begin{equation*}
\begin{CD}
              &  0                        &             &   0                &             &                   &               \\
              &  \downarrow               &             &  \downarrow        &             &                    &               \\
              & \mathcal{O}_{\P^{r+1}}(-1)& =           & \mathcal{O}_{\P^{r+1}}(-1) &             &                     &            \\
              & \downarrow                &             &    \downarrow      &             &                      &               \\
0 \rightarrow &\mathcal{I}_{X/\P^{r+1}}  & \rightarrow & \mathcal{O}_{\P^{r+1}} & \rightarrow & \mathcal{O}_X  & \rightarrow 0 \\
              & \downarrow                &             &  \downarrow        &             & \parallel         &               \\
0 \rightarrow & \mathcal{I}_{X/\P^r}     & \rightarrow & \mathcal{O}_{\P^r} & \rightarrow & \mathcal{O}_{X}    & \rightarrow 0 \\
              & \downarrow                &             &  \downarrow       &             &                    &               \\
              &  0                        &             &   0                 &             &                    & \\
\end{CD}
\end{equation*}
Recall that $H^2 (\P^{r+1},\bigwedge^k \mathcal{M}_{r+1}
(j-1))=H^2 (\P^{r+1},\Omega_{\P^{r+1}} ^k (k+j-1))=0$ for $k+j
\geq 2$. The assertion comes from the cohomology long exact
sequence of the left column.
\end{proof}

\noindent {\bf Remark 2.3.1.} From the exact sequence $0
\rightarrow \mathcal{O}_{\P^{r+1}}(-1) \rightarrow
\mathcal{I}_{X/\P^{r+1}}\rightarrow
\mathcal{I}_{X/\P^r}\rightarrow 0$, it is easily proved that if
the homogeneous ideal of $X \subset \P^{r+e}$ is generated by
forms of degree $\leq d$, then so does $X \subset \P^r$. \qed \\

\begin{corollary}\label{cor:degenerate}
Under the same situation, assume that $X \subset \P^r$ is
nondegenerate.\\
$(1)$ If $H^1 (\P^{r+e},\mathcal{I}_{X/\P^{r+e}} (j)) =0$, then
$H^1 (\P^r,\mathcal{I}_{X/\P^r} (j)) =0$ and hence $X \subset
\P^r$ satisfies $j$-normality.\\
$(2)$ If $H^1 (\P^{r+e} ,\bigwedge^k \mathcal{M}_{r+e} \otimes
\mathcal{I}_{X/\P^{r+e}} (j) )=0$ for $1 \leq k \leq p$ and $j
\geq 2$, then $X \subset \P^r$ satisfies property $N_{2,p}$.
\end{corollary}

\begin{proof}
The assertion comes immediately from Corollary
\ref{cor:criterionN2p} and Lemma \ref{lem:degenrate}.
\end{proof}

\section{The main theorem}
This section is devoted to prove Theorem \ref{thm:veronese}. Let
$X \subset \P (V)=\P^r$ be a nondegenerate projective variety and
let
$\ell \geq 1$ be an integer. Throughout this section we use the following notations:\\

\begin{enumerate}
\item[$\bullet$] $N:={{r+\ell}\choose{\ell}}$ : the dimension of
$S^{\ell} V \cong H^0 (\P^r,\mathcal{O}_{\P^r} (\ell))$
\item[$\bullet$] $Z := \nu_{\ell} (\P (V)) \subset \P^N$
\item[$\bullet$] $V_{\ell}$ : the image of $H^0
(\P^r,\mathcal{O}_{\P^r} (\ell)) \rightarrow H^0 (X,\mathcal{O}_X
(\ell))$\\
\end{enumerate}
Thus we have the following commutative diagram:
\begin{equation*}
\begin{CD}
\quad X              & \quad \quad \subset \quad \quad & \P^r \\
\nu_{\ell} \downarrow     &                      & \downarrow \nu_{\ell} \\
\P (V_{\ell})  & \quad  \subset \quad & \P^N
\end{CD}
\end{equation*}
Let $\mathcal{I}_{X/\P^r}$, $\mathcal{I}_{Z/ \P^N}$,
$\mathcal{I}_{X/\P^N}$ and $\mathcal{I}_{X/ \P (V_\ell)}$ be
sheaves of ideals of $X \subset \P^r$, $Z \subset \P^N$,
$\nu_{\ell}(X) \subset \P^N$ and $\nu_{\ell}(X) \subset \P
(V_{\ell})$, respectively. By snake lemma we have the following:
\begin{equation*}
\begin{CD}
              &                           &             &                    &             & 0                  &               \\
              &                           &             &                    &             & \downarrow         &               \\
              & 0                         &             &                    &             & \mathcal{I}_{X/\P^r} &            \\
              & \downarrow                &             &                    &             & \downarrow         &               \\
0 \rightarrow & \mathcal{I}_{Z / \P^N} & \rightarrow & \mathcal{O}_{\P^N} & \rightarrow & \mathcal{O}_{\P^r} & \rightarrow 0 \\
              & \downarrow                &             & \parallel          &             & \downarrow         &               \\
0 \rightarrow & \mathcal{I}_{X / \P^N}    & \rightarrow & \mathcal{O}_{\P^N} & \rightarrow & \mathcal{O}_{X}    & \rightarrow 0 \\
              & \downarrow                &             &                    &             & \downarrow         &               \\
              & \mathcal{I}_{X/\P^r}    &             &                    &             &    0               &               \\
              & \downarrow                &             &                    &             &                    &               \\
              &   0                       &             &                    &             &                    &               \\
\end{CD}
\end{equation*}
In particular,  we get the short exact sequence
\begin{equation*}
(\lozenge) \quad 0 \rightarrow \mathcal{I}_{Z / \P^N} \rightarrow
\mathcal{I}_{X / \P^N} \rightarrow \mathcal{I}_{X/\P^r}
\rightarrow 0.
\end{equation*}
Now we start to prove our main theorem.\\

\noindent {\bf Proof of Theorem \ref{thm:veronese}.} $(1)$ Since
$Z \subset \P^N$ is projectively normal,
\begin{equation*}
H^1 (\P^N,\mathcal{I}_{Z / \P^N}(k))=0 \quad \mbox{for all $k \geq
1$.}
\end{equation*}
Also from the short exact sequence $0 \rightarrow \mathcal{I}_{Z /
\P^N} \rightarrow \mathcal{O}_{\P^N} \rightarrow
\mathcal{O}_{\P^r} \rightarrow 0$,
\begin{equation*}
H^2 (\P^N,\mathcal{I}_{Z / \P^N}(k))=0\quad \mbox{for all $k \geq
1$.}
\end{equation*}
From the cohomology long exact sequence of $(\lozenge)$, this
implies that
\begin{equation*}
H^1 (\P^N,\mathcal{I}_{X / \P}(k)) \cong H^1 (\P^r,\mathcal{I}_{X
/ \P^r} (k \ell))\quad \mbox{for all $k \geq 1$.}
\end{equation*}
Also our assumption $(\alpha)$ guarantees that $H^1
(\P^r,\mathcal{I}_{X / \P^r} (k \ell))=0$ if $k \geq
\frac{s}{\ell}$. Therefore we have $H^1 (\P^N,\mathcal{I}_{X /
\P^N}(k))=0$ for all $k \geq \frac{s}{\ell}$. This implies that $X
\subset \P (V_{\ell})$ satisfies
$k$-normality for all $k \geq \frac{s}{\ell}$ by Corollary \ref{cor:degenerate}.\\
$(2)$ By Remark 2.2.1 we need to show that the homogeneous ideal
of $\nu_{\ell}(X) \subset \P^N$ is generated by forms of degree
$\leq \mu$ or, equivalently, the multiplication maps
\begin{equation*}
H^0 (\P^N,\mathcal{I}_{X/\P^N} (\mu)) \otimes H^0 (\P^N,
\mathcal{O}_{\P^N} (k)) \rightarrow H^0 (\P^N,\mathcal{I}_{X
/\P^N} (\mu+k))
\end{equation*}
are surjective for all $k \geq 1$ where $\mu =\mbox{max}
\{2,\frac{t}{\ell}\}$. By applying Lemma \ref{lem:mult} to the
short exact sequence $(\lozenge)$ and the line bundle
$\mathcal{O}_{\P^N} (k)$ where $k \geq 1$,
it suffices to check the followings:\\

$(a)$ $\mathcal{I}_{Z/ \P^N} (\mu)$ and $\mathcal{I}_{X/\P^r} (\mu
\ell)$ are globally generated.

$(b)$ The homomorphism
\begin{equation*}
 H^0 (\P^N,\mathcal{I}_{X / \P^N} (k)) \rightarrow H^0
(\P^r,\mathcal{I}_{X/\P^r} (k \ell))
\end{equation*}

is surjective for all $k \geq \mu$.

$(c)$ The multiplication maps $$\begin{cases} (\triangleleft)
\quad H^0 (\P^N,\mathcal{I}_{Z / \P^N} (\mu)) \otimes H^0
(\P^N,\mathcal{O}_{\P^N} (k)) \rightarrow H^0
(\P^N,\mathcal{I}_{Z / \P^N} (\mu+k))\\
(\triangleright) \quad H^0 (\P^r,\mathcal{I}_{X/\P^r}(\mu \ell))
\otimes H^0 (\P^N,\mathcal{O}_{\P^N} (k)) \rightarrow H^0 (\P^r
,\mathcal{I}_{X/\P^r} ((\mu+k)\ell))
\end{cases}$$

are surjective for all $k \geq 1$.\\

\noindent Since the homogeneous ideal of $Z \subset \P^N$ is
generated by quadrics, $\mathcal{I}_{Z/ \P^N} (\mu)$ is globally
generated. Also the homomorphism in $(b)$ and the multiplication
map $(\triangleleft)$ are surjective for all $k \geq 1$ since $Z
\subset \P^N$ is projectively normal. The global generation of
$\mathcal{I}_{X/\P^r} (\mu \ell)$ and the surjectivity of the
multiplication map $(\triangleright)$ come from the assumption
$(\beta)$. \\
$(3)$ Let $\mathcal{M}_{\ell}$ be the restriction of
$\mathcal{M}_{\P^N} := \Omega_{\P^N} (1)$ to $Z$. Clearly
$\mathcal{M}_{\ell}$ is the kernel of the evaluation homomorphism
$V_{\ell} \otimes \mathcal{O}_{\P^r} \rightarrow
\mathcal{O}_{\P^r} (\ell) \rightarrow 0$ and the short exact
sequence
\begin{equation*}
0 \rightarrow \mathcal{M}_{\ell} \rightarrow V_{\ell} \otimes
\mathcal{O}_{\P^r} \rightarrow \mathcal{O}_{\P^r} (\ell)
\rightarrow 0
\end{equation*}
is the restriction of the Euler sequence on $\P^N$ to  $Z$. From
the short exact sequence $(\lozenge)$, we have the long exact
sequence
\begin{equation*}
H^1 (\P^N,\bigwedge^k \mathcal{M}_{\P^N} \otimes \mathcal{I}_{Z /
\P^N} (j)) \rightarrow H^1 (\P^N,\bigwedge^k \mathcal{M}_{\P^N}
\otimes \mathcal{I}_{X / \P^N} (j)) \rightarrow H^1
(\P^r,\bigwedge^k \mathcal{M}_{\ell} \otimes \mathcal{I}_{X /
\P^r} (j\ell))
\end{equation*}
of cohomology groups. Since $Z \subset \P^N$ satisfies property
$N_{\ell}$ by Theorem \ref{thm:GreenVeronese},
\begin{equation*}
H^1 (\P^N,\bigwedge^k \mathcal{M}_{\P^N} \otimes \mathcal{I}_{\P^r
/ \P^N} (j))=0
\end{equation*}
for $1 \leq k \leq \ell$ and all $j \geq 2$ by Corollary
\ref{cor:criterionN2p}. Also we claim that
\begin{equation*}
H^1 (\P^r,\bigwedge^k \mathcal{M}_{\ell} \otimes \mathcal{I}_{X /
\P^r} (j\ell))=0
\end{equation*}
for $1 \leq k \leq \ell +1$ and $j\ell \geq 2\ell \geq m+k-1$.
Indeed $\bigwedge^k \mathcal{M}_{\ell}$ is $(\ell+1)$-regular with
respect to $\mathcal{O}_{\P^r} (1)$ by the next Proposition
\ref{prop:regularityVeronese} since $(\P^r,\mathcal{O}_{\P^r}
(\ell))$ satisfies property $N_{\ell}$ by Green's theorem.
Therefore $\bigwedge^k \mathcal{M}_{\ell}\otimes \mathcal{I}_{X /
\P^r}$ is $(m+\ell+1)$-regular with respect to $\mathcal{O}_{\P^r}
(1)$. For details, we refer the reader to see Proposition 1.8.9 in
\cite{L2}. This completes the proof of our claim.\\

\underline{\textit{Case 1.}} Assume that $\frac{m+1}{2} \leq \ell
< m$. Then $2\ell \geq m+k-1$ holds for $1 \leq k \leq 2\ell-m+1$.
Therefore
\begin{equation*}
H^1 (\P^N,\bigwedge^k \mathcal{M}_{\P^N} \otimes \mathcal{I}_{X /
\P^N} (j))=0 \quad \mbox{for $1 \leq k \leq 2\ell-m+1$ and all $j
\geq 2$.}
\end{equation*}
Consequently Corollary \ref{cor:degenerate} shows
that $X \subset \P (V_{\ell})$ satisfies property $N_{2,2\ell-m}$.\\

\underline{\textit{Case 2.}} Assume that $\ell \geq m$. Then
$2\ell \geq m+k-1$ holds for $1 \leq k \leq \ell+1$ and hence we
have
\begin{equation*}
H^1 (\P^N,\bigwedge^k \mathcal{M}_{\P^N} \otimes \mathcal{I}_{X /
\P^N} (j))=0 \quad \mbox{for $1 \leq k \leq \ell$ and $j \geq 2$.}
\end{equation*}
Therefore $X (X) \subset \P (V_{\ell})$ satisfies property
$N_{2,\ell}$ by Corollary \ref{cor:degenerate}. \qed \\

\begin{proposition}\label{prop:regularityVeronese}
If $(\P^r,\mathcal{O}_{\P^r} (\ell))$ satisfies property $N_p$,
then $\bigwedge^k \mathcal{M}_{\ell}$ is $(\ell+1)$-regular for
every $1 \leq k \leq p+1$.
\end{proposition}

\begin{proof}
By definition of the Castelnuovo-Mumford regularity, we need to
show that
\begin{equation*}
H^i (\P^r,\bigwedge^k \mathcal{M}_{\ell} (\ell+1-i))=0
\end{equation*}
for $1 \leq k \leq p+1$ and all $i \geq 1$. For each $k \geq 0$,
consider the short exact sequence
\begin{equation*}
(\star) \quad 0 \rightarrow \bigwedge^{k+1} \mathcal{M}_{\ell}
\rightarrow \bigwedge^{k+1} V_{\ell} \otimes \mathcal{O}_{\P^r}
\rightarrow \bigwedge^k \mathcal{M}_{\ell} (\ell) \rightarrow 0.
\end{equation*}

\underline{\textit{Step 1.}} Since $H^r (\P^r,\mathcal{O}_{\P^r}
(j))=0$ for all $j \geq -r$, we have
\begin{equation*}
H^r (\P^r,\bigwedge^k \mathcal{M}_{\ell} (j))=0
\end{equation*}
for all $k \geq 1$ and $j \geq \ell-r$ by the cohomology long
exact sequence induced from $(\star)$.\\

\underline{\textit{Step 2.}} If $(\P^r,\mathcal{O}_{\P^r} (\ell))$
satisfies property $N_p$, then
\begin{equation*}
(\star\star) \quad H^1 (\P^r,\bigwedge^k \mathcal{M}_{\ell}
(\ell))  =0
\end{equation*}
for every $1 \leq k \leq p+1$ by the cohomological criterion in
\cite{EL}.\\

\underline{\textit{Step 3.}} By $(\star)$ and $(\star\star)$, we
have
\begin{equation*}
H^1 (\P^r,\bigwedge^k \mathcal{M}_{\ell} (\ell)) \cong H^2
(\P^r,\bigwedge^{k+1} \mathcal{M}_{\ell}) =0
\end{equation*}
for every $1 \leq k \leq p+1$. Then Proposition 1.7 in \cite{OP}
guarantees that
\begin{equation*}
(\star\star\star) \quad H^2 (\P^r,\bigwedge^{k+1}
\mathcal{M}_{\ell} (t)) \cong H^1 (\P^r,\bigwedge^k
\mathcal{M}_{\ell} (\ell+t)) =0
\end{equation*}
for every $1 \leq k \leq p+1$ and $t \geq 0$.\\

\underline{\textit{Step 4.}} Assume that $1 \leq i \leq r-1$. If
$k \geq i$, then
\begin{equation*} H^i (\P^r,\bigwedge^k \mathcal{M}_{\ell} (\ell+1-i))
\cong \cdots \cong H^1 (\P^r,\bigwedge^{i-k+1} \mathcal{M}_{\ell}
(\ell+(\ell-1)(i-1)))=0
\end{equation*}
by $(\star\star\star)$. Also if $k < i$, then
\begin{equation*}
H^i (\P^r,\bigwedge^k \mathcal{M}_{\ell} (\ell+1-i)) \cong \cdots
\cong H^{i-k} (\P^r,\mathcal{O}_{\P^r} (\ell+1-i+k\ell))=0.
\end{equation*}\\
Therefore $\bigwedge^k \mathcal{M}_d$ is $(\ell+1)$-regular for
every $1 \leq k \leq p+1$.
\end{proof}

\noindent {\bf Proof of Theorem \ref{thm:main1}.} It is well known
that if $X \subset \P H^0 (X,L)$ is $m$-regular, then
$k$-normality holds for all $k \geq m-1$ and the homogeneous ideal
is generated by forms of degree $\leq m$ (Lecture 14 in \cite{M}).  \\
$(1)$ Let $V_{\ell}$ be the image of $H^0 (\P^r,\mathcal{O}_{\P^r}
(\ell)) \rightarrow H^0 (X,L^{\ell})$. If $\ell \geq
\frac{m-1}{2}$, then Theorem \ref{thm:veronese}.$(1)$ implies that
$X \subset \P (V_{\ell})$ is $k$-normal for all $k \geq 2$.
Therefore the linearly normal embedding $X \subset \P H^0
(X,L^{\ell})$ is projectively normal.\\
$(2)$ Since $(m-1)$-normality holds for $X \subset \P H^0 (X,L)$,
$V_{m-1}=H^0 (X,L^{m-1})$. Also Theorem \ref{thm:veronese}.$(3)$
implies that $X \subset \P H^0 (X,L^{m-1})$ satisfies property
$N_{2,m-2}$. Therefore
$(X,L^{m-1})$ satisfies property $N_{m-2}$.\\
$(3)$ Since $X \subset \P H^0 (X,L)$ is $\ell$-normal for all
$\ell \geq m$, $V_{\ell}=H^0 (X,L^{\ell})$. Also Theorem
\ref{thm:veronese}.$(3)$ implies that $X \subset \P H^0
(X,L^{\ell})$ satisfies property $N_{2,\ell}$ if $\ell \geq m$.
Therefore $(X,L^{\ell})$ satisfies property $N_{\ell}$. \qed \\

\noindent {\bf Remark 3.1.} In many cases, an $m$-regular
projective variety $X \subset \P H^0 (X,L)$ satisfies
$s$-normality for some $s < m-1$. Suppose that $s$-normality holds
for some $\frac{m+1}{2} \leq s \leq m-2$. Then by the same reason
as in the above proof, Theorem \ref{thm:veronese}.$(3)$ guarantees
that $(X,L^s)$ satisfies
property $N_{2s-m}$. \qed \\

\section{Applications and examples}
In this section, we present some applications of our results.
\subsection{Generalized Ottaviani-Paoletti's conjecture} By
M. Green's Theorem 2.2 in \cite{Green2},
$(\P^r,\mathcal{O}_{\P^r}(\ell))$ satisfies property $N_{\ell}$.
We generalize this result to arbitrary projective varieties in
Theorem \ref{thm:main1}. On the other hand, G. Ottaviani and R.
Paoletti\cite{OP} prove that if $r \geq 2$ and $\ell \geq 3$, then
$(\P^r,\mathcal{O}_{\P^r}(\ell))$ fails to satisfy property
$N_{3\ell-2}$. Also they conjectured that for $r \geq 2$ and $\ell
\geq 3$, $(\P^r,\mathcal{O}_{\P^r}(\ell))$ satisfies property
$N_p$ if and only if $p \leq 3\ell-3$. The cases $r \geq 3$ are
still open except $r=\ell=3$. We generalize Ottaviani-Paoletti's
conjecture to very ample line bundles on an arbitrary projective
variety. See ``Conjecture A" and ``Conjecture B" in $\S 1$. And
the aim of this subsection is to prove that if Conjecture A is
true, then Conjecture B is also true. \\

\noindent {\bf Proof of Theorem \ref{thm:imply}.} We use the short
exact sequence $(\lozenge)$ in $\S 3$. Note that by Corollary
\ref{cor:degenerate} we need to show that
\begin{equation*}
H^1 (\P^N,\bigwedge^k \mathcal{M}_{\P^N} \otimes
\mathcal{I}_{X/\P^N} (j) )=0 \quad \mbox{for $1 \leq k \leq p$ and
$j \geq 2$.}
\end{equation*}
The short exact sequence $(\lozenge)$ gives the cohomology long
exact sequence
\begin{equation*}
H^1 (\P^N,\bigwedge^k \mathcal{M}_{\P^N} \otimes \mathcal{I}_{Z /
\P^N} (j)) \rightarrow H^1 (\P^N,\bigwedge^k \mathcal{M}_{\P^N}
\otimes \mathcal{I}_{X / \P^N} (j)) \rightarrow H^1
(\P^r,\bigwedge^k \mathcal{M}_{\ell} \otimes \mathcal{I}_{X /
\P^r} (j\ell)).
\end{equation*}
Since $(\P^r,\mathcal{O}_{\P^r}(\ell))$ satisfies property $N_p$,
\begin{equation*}
H^1 (\P^N,\bigwedge^k \mathcal{M}_{\P^N} \otimes \mathcal{I}_{Z /
\P^N} (j))=0 \quad \mbox{for $1 \leq k \leq p$ and $j \geq 2$.}
\end{equation*}
Also $\bigwedge^k \mathcal{M}_{\ell}$ is $(\ell+1)$-regular for
every $1 \leq k \leq p+1$ by Proposition
\ref{prop:regularityVeronese}. Since $\mathcal{I}_{X/\P^r}$ is
$m$-regular, $\bigwedge^k \mathcal{M}_{\ell} \otimes
\mathcal{I}_{X/\P^r}$ is $(m+\ell+1)$-regular by Proposition 1.8.9
in \cite{L2}. Therefore $H^1 (\P^r,\bigwedge^k \mathcal{M}_{\ell}
\otimes \mathcal{I}_{X/\P^r} (j\ell))=0$ for $1 \leq k \leq p$ and
$j\geq2$ since we assume that $\ell \geq m$. As a consequence, we
have
\begin{equation*}
H^1 (\P^N,\bigwedge^k \mathcal{M}_{\P^N} \otimes \mathcal{I}_{X /
\P^N} (j))=0 \quad \mbox{for $1 \leq k \leq p$ and $j \geq 2$.}
\end{equation*}
Therefore $X (X) \subset \P (V_{\ell})$ satisfies property
$N_{2,\ell}$ by Corollary \ref{cor:degenerate}. \qed \\

\subsection{Castelnuovo-Mumford regularity of projective varieties}
Our results are stated in terms of the Castelnuovo-Mumford
regularity of projective varieties. Here we shortly review known
results in this field. We refer to \cite{K1},\cite{K2} and
\cite{K3} for details. There is a well known conjecture due to
Eisenbud and Goto\cite{EG} which gives a bound for regularity in
terms of the degree and the codimension:\\

{\bf Regularity conjecture.} Let $X \subset \P^r$ be a
nondegenerate integral projective variety of degree $d$ and
codimension $e$. Then it satisfies $(d-e+1)$-regularity.\\

\noindent Let us recall the following works developed toward this conjecture.\\

\begin{enumerate}
\item[$\bullet$] (L. Gruson, R. Lazarsfeld and C. Peskine,
\cite{GLP}) Let $C \subset \P^r$ be a nondegenerate integral curve
of degree $d$. Then $C$ is $(d+2-r)$-regular. \item[$\bullet$] (R.
Lazarsfeld, \cite{L1}) Let $S \subset \P^r$ be a nondegenerate
smooth surface of degree $d$. Then $S$ is $(d+3-r)$-regular.
\item[$\bullet$] (Ziv Ran, \cite{Ran}) Let $X \subset \P^r$,
$r\geq 9$,  be a nondegenerate smooth threefold of degree $d$.
Then $X$ is $(d+4-r)$-regular. \item[$\bullet$] (Sijong Kwak,
\cite{K1}\cite{K2}\cite{K3}) Let $X \subset \P^r$ be a
nondegenerate smooth projective variety of degree $d$.
\begin{enumerate}
\item[$1.$] If $X$ is a threefold and $r=5$, then $X$ is
$(d-1)$-regular. \item[$2.$] If $X$ is a threefold and $6 \leq r
\leq 8$, then $X$ is $(d+5-r)$-regular. \item[$3.$] If $X$ is a
fourfold, then $X$ is $(d+9-r)$-regular. \item[$4.$] If $X$ is a
fivefold, then $X$ is $(d+16-r)$-regular. \item[$5.$] If $X$ is a
sixfold, then $X$ is $(d+27-r)$-regular.
\end{enumerate}
\item[$\bullet$] (A. Bertram, L. Ein and R. Lazarsfeld,
\cite{BEL}) Let $X \subset \P^r$ be a nondegenerate smooth
projective variety. Let $n$ be the dimension of $X$ and let $d$ be
the degree of $X \subset \P^r$. Then $X$ is
$(\mbox{min}\{r-n,n+1\}(d-1)-n+1)$-regular.\\
\end{enumerate}

\noindent Note that Theorem \ref{thm:main1} is stated in terms of
the Castelnuovo-Mumford regularity of the linearly normal
embedding. For a projective variety $X$ of dimension $n$ and a
very ample line bundle $L \in \mbox{Pic}X$ of degree $d$, the
$\Delta$-genus $\Delta(X,L)$ of $(X,L)$ is defined to be $d+n-h^0
(X,L)$. Therefore for linearly
normal varieties, the Eisenbud-Goto's regularity conjecture is\\

{\bf Regularity conjecture for linearly normal varieties.} Let $X$
be an integral projective variety and let $L$ be a very ample line
bundle on $X$. Then the linearly normal embedding $X \subset
\P H^0 (X,L)$ satisfies $(\Delta(X,L)+2)$-regularity.\\

One can apply Theorem \ref{thm:main1} to the above works about
Castelnuovo-Mumford regularity. For an example, we obtain the
following

\begin{corollary}\label{cor:surfaces}
Let $S$ be a smooth surface and let $L \in \mbox{Pic}S$ be a very
ample line bundle. Put $m = \Delta(S,L)+2$. Then \\
$(1)$ $L^{\ell}$ is normally generated for $\ell \geq
\frac{m-1}{2}$.\\
$(2)$ $L^{m-1}$ satisfies property $N_{m-2}$.\\
$(3)$ $L^{\ell}$ satisfies property $N_{\ell}$ for $\ell \geq m$.
\end{corollary}

\begin{proof}
By R. Lazarsfeld's result in \cite{L1}, $S \subset \P H^0 (S,L)$
is $(\Delta(S,L)+2)$-regular. Thus the assertion follows
immediately from Theorem \ref{thm:main1}.
\end{proof}

\subsection{$3$-regular varieties} Let $L$ be a very ample line
bundle on a smooth projective variety $X$ where the linearly
normal embedding $X \subset \P H^0 (X,L)$ is $3$-regular. Then it
is projectively normal and Theorem \ref{thm:main1} implies that
 $$ (\dag) \begin{cases} (i)& (X,L)\mbox{ satisfies property
 $N_0$.}\\
(ii) & (X,L^2)\mbox{ satisfies property $N_1$.} \\
(iii) &(X,L^3)\mbox{ satisfies property $N_4$}.\\
(iv) &(X,L^{\ell}),~\ell \geq 4, \mbox{ satisfies property
$N_{\ell}$}.\end{cases}$$ For $(iii)$, see Example 1.1. Here we
exhibits examples where $3$-regularity holds.  \\

$(4.3.1)$ Let $S$ be a K3 surface. Then any linearly normal
embedding of $S$ is $3$-regular. Indeed for a very ample line
bundle $L$ on $S$, let $C \in  |L|$ be a smooth section. Since
$H^1 (S,L^j)=0$ for all $j \geq 0$, $S \subset \P H^0 (S,L)$ and
$C \subset \P H^0 (C,L|_C)$ has the same Betti table. Thus
Noether's theorem implies that $S \subset \P H^0 (S,L)$ is
projectively normal and hence it is $3$-regular. Therefore $(S,L)$
satisfies $(\dag)$. We should mention Gallego-Purnaprajna's result
in \cite{GP2}. Recall that $L^2 \geq 4$, $L^2=4$ if and only if
$(S,L)$ defines a quartic hypersurface in $\P^3$, and $L^2=6$ if
and only if $(S,L)$ defines the complete intersection of a quadric
equation and a cubic equation in $\P^4$. Therefore $L^2 \geq 8$
except these two cases.

\begin{theorem}[F. J. Gallego and B. P. Purnaprajna]
Let $S$ be a K3 surface and let $L \in \mbox{Pic}S$ be a very
ample line bundle.\\
$(1)$ If $L^2=4$ or $6$, then $(S,L^{\ell})$ satisfies property
$N_{\ell -1}$.\\
$(2)$ If $L^2 \geq 8$ and the general member of $|L|$ is
non-trigonal, then $(S,L^{\ell})$ satisfies property $N_{\ell}$.
\end{theorem}

Indeed their result is about the case when $L$ is ample and base
point free. Clearly $(\dag)$ refines or reproves their result for $\ell \geq 3$.\\

$(4.3.2)$ Let $S \subset \P^{g-1}$ be a smooth linearly normal
Enriques surface. Then $g \geq 6$. $\mathcal{O}_S (1)$ is called a
Reye polarization if $g=6$ and $S$ lies on a quadric. In
\cite{GLM}, the projective normality of Enriques surfaces is
studied.

\begin{theorem}[L. Giraldo, A. F. Lopez and R.
Mu$\widetilde{n}$oz]\label{thm:Enriques} $(1)$ If $g=6$ and
$\mathcal{O}_S (1)$ is a Reye polarization, then it is $j$-normal
for every $j \geq 3$, $4$-regular, and its
homogeneous ideal is generated by forms of degree $\leq 3$.\\
$(2)$ If $g \geq 7$ or $g=6$ and $\mathcal{O}_S (1)$ is not a Reye
polarization, then $S \subset \P^{g-1}$ is $3$-regular.
\end{theorem}

Therefore we have

\begin{corollary}\label{cor:Enriques}
Let $S \subset \P^{g-1}$ be a smooth linearly normal Enriques
surface.\\
$(1)$ If $g=6$ and $\mathcal{O}_S (1)$ is a Reye polarization,
then $\mathcal{O}_S (2)$ satisfies property $N_1$, $\mathcal{O}_S
(3)$ satisfies property $N_2$, and $\mathcal{O}_S (\ell)$, $\ell
\geq 4$,
satisfies property $N_{\ell}$.\\
$(2)$ If $g \geq 7$ or $g=6$ and $\mathcal{O}_S (1)$ is not a Reye
polarization, then $\mathcal{O}_S (1)$ satisfies property $N_0$,
$\mathcal{O}_S (2)$ satisfies property $N_1$, $\mathcal{O}_S (3)$
satisfies property $N_4$, and $\mathcal{O}_S (\ell)$, $\ell \geq
4$, satisfies property $N_{\ell}$.
\end{corollary}

\begin{proof}
$(1)$ By the Riemann-Roch theorem, the degree of $S \subset
\P^{g-1}$ is $10$. Therefore Theorem 2.14 in \cite{GP} guarantees
that $\mathcal{O}_S (2)$ is satisfies property $N_1$ and
$\mathcal{O}_S (3)$ satisfies property $N_2$. Also $\mathcal{O}_S
(\ell)$, $\ell \geq 4$, satisfies property $N_{\ell}$ by Theorem
\ref{thm:main1}.\\
$(2)$ This follows immediately from Theorem \ref{thm:Enriques} and
$(\dag)$.
\end{proof}

\subsection{Complex torus} In \cite{Rubei1} and \cite{Ru1}, E. Rubei studies syzygies of a power of
normally generated very ample line bundles on a complex torus. She
prove that $(1)$ for a normally generated line bundle $L \in
\mbox{Pic}X$ on a complex torus $X$, $(X,L^{\ell})$ satisfies
property $N_{\ell-1}$, and $(2)$ for a normally presented line
bundle $L \in \mbox{Pic}X$ on $X$, $(X,L^{\ell})$ satisfies
property $N_{\ell}$. This result is refined by Theorem
\ref{thm:main1} as follows:

\begin{corollary}\label{cor:refineRubei}
Let $X$ be a complex torus of dimension $n$ and let $L \in
\mbox{Pic}X$ be a very ample line bundle such that $X \subset \P
H^0
(X,L)$ satisfies $(n+1)$-normality. Then\\
$(1)$ $L^{n+1}$ satisfies property $N_n$.\\
$(2)$ $L^{\ell}$ satisfies property $N_{\ell}$ for $\ell \geq
n+2$.
\end{corollary}

\begin{proof}
Since $(n+1)$-normality holds for $X \subset \P H^0 (X,L)$, it is
$(n+2)$-regular by Kodaira vanishing theorem. Thus Theorem
\ref{thm:main1} guarantees the assertion.
\end{proof}

\subsection{A remark on a Rubei's result} In \cite{Rubei}, Rubei
consider the problem to see how property $N_p$ propagates through
powers. More precisely she prove that if $(X,L)$ satisfies
property $N_p$, then $(X,L^{\ell})$ satisfies property
$N_{\mbox{min}\{\ell,p\}}$. In this subsection, we prove the
following:

\begin{theorem}\label{thm:Rubeitype}
Let $X$ be a projective variety and let $L \in \mbox{Pic}X$ be a
very ample line bundle such that $(X,L)$ satisfies property $N_p$.
Then for $\ell \geq 2$, $(X,L^{\ell})$ satisfies property $N_k$
where $$k= \begin{cases} \mbox{min}\{5,p\} & \mbox{for $\ell=2$,}
\\ \mbox{min}\{4,p\} & \mbox{for $\ell=3$, and} \\ \mbox{min}\{\ell,p\} & \mbox{for $\ell \geq 4$.} \end{cases}$$
\end{theorem}

\begin{proof}
We use the short exact sequence $(\lozenge)$ in $\S 3$. As in the
proof of Theorem \ref{thm:veronese}.$(3)$, we have the long exact
sequence
\begin{equation*}
H^1 (\P^N,\bigwedge^k \mathcal{M}_{\P^N} \otimes \mathcal{I}_{Z /
\P^N} (j)) \rightarrow H^1 (\P^N,\bigwedge^k \mathcal{M}_{\P^N}
\otimes \mathcal{I}_{X / \P^N} (j)) \rightarrow H^1
(\P^r,\bigwedge^k \mathcal{M}_{\ell} \otimes \mathcal{I}_{X /
\P^r} (j\ell))
\end{equation*}
of cohomology groups. Since $(X,L)$ satisfies property $N_p$, we
have
\begin{equation*}
H^1 (\P^r,\bigwedge^k \mathcal{M}_{\ell} \otimes \mathcal{I}_{X /
\P^r} (j\ell))=0 \quad \mbox{for $1 \leq k \leq p$ and $j \ell
\geq 2$.}
\end{equation*}
For the desired vanishing of $H^1 (\P^N,\bigwedge^k
\mathcal{M}_{\P^N} \otimes \mathcal{I}_{Z / \P^N} (j))$, see
Theorem \ref{thm:GreenVeronese}, Example $1.1$ and Example $1.2$.
Also we refer the reader to see \cite{OP}.
\end{proof}

\noindent {\bf Remark 4.5.1.} As in the case of Theorem
\ref{thm:main1}, if Ottaviani-Paoletti's conjecture is true then
Theorem \ref{thm:Rubeitype} is automatically refined. That is, for
$\ell \geq 3$, $(X,L^{\ell})$ satisfies property
$N_{\mbox{min}\{3\ell-3,p\}}$. \qed \\

\noindent {\bf Remark 4.5.2.} As remarked in \cite{Rubei} and
\cite{Rubei2} it is not true that if $(X,L)$ satisfies property
$N_p$ then any power of $L$ satisfies property $N_p$. Indeed Rubei
exhibits a concrete example. Also one can find more examples in
\cite{P}. More precisely, let $(X,L)$ be a ruled scroll over a
smooth curve $C$ with the projection morphism $\pi : X \rightarrow
C$. For an arbitrary $p \geq 1$, one can find $(X,L)$ which
satisfies property $N_p$. Now assume that $\mbox{dim}(X) \geq 3$.
Then for $\ell \geq 3$, $(X,L^{\ell})$ fails to satisfy property
$N_{3\ell-2}$. Therefore if $p \geq 3\ell-2$, then $(X,L)$
satisfies property $N_p$ while $(X,L^{\ell})$ fails to satisfy
property $N_p$. For details, see Corollary 3.7 in \cite{P}. \qed \\

\subsection{Complete intersection} Let $X \subset \P^r$ be a normal variety of dimension $\geq 1$
which is a complete intersection. Then it is projectively normal.
For the property $N_1$ of $\mathcal{O}_X (\ell)$, we have the
following

\begin{corollary}\label{cor:ci}
Under the situation just stated, assume that the homogeneous ideal
of $X \subset \P^r$ is generated by forms $F_1,\cdots,F_e$ of
degrees $d_1 \leq d_2 \leq \cdots \leq d_e$ where
$e=\mbox{codim}(X,\P^r)$. Then $(X,\mathcal{O}_X (\ell))$
satisfies property $N_1$ if and only if $\ell \geq \frac{d_e}{2}$.
\end{corollary}

\begin{proof}
If $\ell \geq \frac{d_e}{2}$, then $\mathcal{O}_X (\ell)$
satisfies property $N_1$ by Theorem \ref{thm:veronese}.
Conversely, assume that $\mathcal{O}_X (\ell)$ satisfies property
$N_1$. We follow notations in $\S 3$. From the short exact
sequence
\begin{equation*}
(\lozenge) \quad 0 \rightarrow \mathcal{I}_{Z / \P^N} \rightarrow
\mathcal{I}_{X / \P^N} \rightarrow \mathcal{I}_{X/\P^r}
\rightarrow 0.
\end{equation*}
if $\mathcal{I}_{X / \P^N} (2)$ is globally generated, then
$\mathcal{I}_{X/\P^r} (2\ell)$ is also globally generated. Since
$\mathcal{I}_{X/\P^r} (2\ell)$ is globally generated if and only
if $2\ell \geq d_e$, the converse is proved.
\end{proof}

\end{document}